\documentclass[12pt]{amsart}

\usepackage{amssymb}

\usepackage{enumitem}

\usepackage{graphicx}

\makeatletter
\@namedef{subjclassname@2020}{%
  \textup{2020} Mathematics Subject Classification}
\makeatother

\usepackage[T1]{fontenc}
\usepackage{tikz-cd}
\usepackage[nameinlink]{cleveref}

\newtheorem{thm}{Theorem}[section]
\newtheorem{cor}[thm]{Corollary}

\newtheorem{prop}[thm]{Proposition}
\newtheorem{prob}[thm]{Problem}

\theoremstyle{definition}
\newtheorem{Def}[thm]{Definition}
\newtheorem{obs}[thm]{Remark}
\newtheorem{ex}[thm]{Example}

\numberwithin{equation}{section}

\newcommand{\N}{\mathbb{N}}
\newcommand{\Z}{\mathbb{Z}}
\newcommand{\charc}{\: \raisebox{2pt}{$\chi$}}
\newcommand{\R}{\mathbb{R}}
\newcommand{\B}{\mathbb{B}}
\newcommand{\LF}{\mathcal{F}}
\newcommand{\BL}{\mathcal{L}}

\newcommand{\norm}[1]{\left\| #1 \right\|}
\newcommand{\hprod}[1]{\left\langle #1 \right\rangle}

\DeclareMathOperator{\spanset}{span}

\DeclareMathOperator{\fin}{fin}

\newcommand{\cspan}{\overline{\spanset}}
\DeclareMathOperator{\Lip}{Lip}

\newcommand{\isometric}{\cong}

\makeatletter
\newcommand*{\rest}[2]{%
   \ensuremath{%
        \mathpalette{\@resttwo{#1}}{#2}%
   }%
 }   
 \newcommand*{\@resttwo}[3]{%
   \begingroup
     #2%
     \sbox0{$\m@th#2\vphantom{\big|}#1$}%
     \setbox4\vbox{\hbox{$\m@th#2\upharpoonright$}\kern\z@}%
     \setbox6\hbox{$#2\vcenter{}$}%
     \sbox4{\kern-.4\wd4\box4}%
     \dimen0=\ht0 %
     \advance\dimen0 by -\ht6 %
     \dimen2=\dp0 %
     \advance\dimen2 by \ht6 %
     \ifdim\dimen2>\dimen0 %  
       \dimen0=\dimen2 %
     \else
       \dimen0=\dimen0 %
     \fi
     \dimen2=\ht6 %
     \advance\dimen2 by -\dimen0 %
     \dimen0=2\dimen0 %
     \def\DelimCorr{%  
       \mskip.5\thinmuskip
       \nonscript\mskip.5\thinmuskip
     }%
     \begingroup
       #2#1%
     \endgroup
     \mathclose{%
       \DelimCorr
       \raisebox{\dimen2}{\resizebox{!}{\dimen0}{\box4}}%
        }_{#3}%
   \endgroup
 }
\makeatother

\makeatletter
\newcommand*{\ext}[2]{%
   \ensuremath{%
        \mathpalette{\@exttwo{#1}}{#2}%
   }%
 }   
 \newcommand*{\@exttwo}[3]{%
   \begingroup
     #2%
     \sbox0{$\m@th#2\vphantom{\big|}#1$}%
     \setbox4\vbox{\hbox{$\m@th#2\downharpoonright$}\kern\z@}%
     \setbox6\hbox{$#2\vcenter{}$}%
     \sbox4{\kern-.4\wd4\box4}%
     \dimen0=\ht0 %
     \advance\dimen0 by -\ht6 %
     \dimen2=\dp0 %
     \advance\dimen2 by \ht6 %
     \ifdim\dimen2>\dimen0 %  
       \dimen0=\dimen2 %
     \else
       \dimen0=\dimen0 %
     \fi
     \dimen2=\ht6 %
     \advance\dimen2 by -\dimen0 %
     \dimen0=2\dimen0 %
     \def\DelimCorr{%  
       \mskip.5\thinmuskip
       \nonscript\mskip.5\thinmuskip
     }%
     \begingroup
       #2#1%
     \endgroup
     \mathclose{%
       \DelimCorr
       \raisebox{\dimen2}{\resizebox{!}{\dimen0}{\box4}}%
        }^{#3}%
   \endgroup
 }
\makeatother
\begin{document}

%%%%% To ease editing, for IMPAN journals add:

\baselineskip=17pt

%%%%%%%%%%%%%%%%

\title[Borel Measures which Induce Lipschitz free Space Elements]{A Characterization of Borel Measures which Induce Lipschitz free Space Elements}

\author{Lucas Maciel Raad}
\address{Instituto de Ci\^encia e Tecnologia da Universidade federal de S\~ao Paulo, Av. Cesare Giulio Lattes, 1201, ZIP 12247-014 S\~ao Jos\'e dos Campos/SP, Brasil}
\email{raad.lucas@unifesp.br}

\date{}

\begin{abstract}
    We solve a problem of Aliaga and Pernecká about Lipschitz free spaces (denoted by $\LF(M)$):
    \begin{quote}
        Does every Borel measure $\mu$ on a complete metric space $M$ such that $\int d(m,0) d |\mu|(m)< \infty$ induce a weak$^*$ continuous functional $\BL\mu \in \LF(M)$ by the mapping $\BL\mu(f)=\int f d \mu$?
    \end{quote}
    In particular, we obtain a characterization of the Borel measures $\mu$ such that $\BL\mu \in \LF(M)$, which indeed implies inner-regularity for complete metric spaces. We also prove that every Borel measure on $M$ induces an element of $\LF(M)$ if and only if the weight of $M$ is strictly less than the least real-valued measurable cardinal, and thus the existence of a metric space on which there is a measure $\mu$ such that ${\BL\mu \in \LF(M)^{**} \setminus \LF(M)}$ cannot be proven in ZFC. Finally, we partially solve a problem of Aliaga on whether every sequentially normal functional on $\Lip_0(M)$ is normal.
\end{abstract}



\maketitle

\section{Introduction}

    The Lipschitz free space over a metric space $M$, denoted by $\LF(M)$, is a Banach space that encodes in its linear structure the metric structure of $M$.  It is related to many areas of mathematics, including metric geometry (\cite{Gartland}), harmonic analysis (\cite{Doucha_2022}) and machine learning (\cite{Luxburg}). We will focus on its relation to measure theory, exemplified by its construction via the closed span of Dirac Measures (\cite{GodefroyKalton}) and by its relation to the mass transfer problem (\cite{Weaver} Chapter 3.3). More specifically we will focus on a problem in  \cite{Aliaga_integral} regarding a theorem similar to the Riesz–Markov–Kakutani Representation Theorem.
    
    From \cite{Aliaga_integral} Proposition 4.3\footnote{The original proposition refers to the more general case of measures on a compactification of $M$ named the uniform compactification, we will not need this level of generality. } we know that if $\mu$ is a Borel measure on a pointed complete metric space $M$, then the mapping $\BL\mu$ defined by $\BL\mu(f)=\int f d \mu$ is an element of $\LF(M)^{**}$ if and only if $\int d(m,0) d |\mu|(m)< \infty$; and from \cite{Aliaga_integral} Proposition 4.4 we know that if $\mu$ is also inner-regular or $M$ is separable then $\BL\mu \in \LF(M)$. Problem 2 in \cite{Aliaga_integral} asks if the hypothesis that $\mu$ is inner-regular is indeed necessary for $\BL\mu$ to be an element of $\LF(M)$.
    
    We will show that if $M$ is a metric space and $\mu$ is a Borel measure on $M$, then we have that $\BL\mu \in \LF(M)$ if and only if $\int d(m,0) d |\mu|(m)< \infty$ and $\mu$ is concentrated on a separable subset of $M$. Which, when $M$ is complete, indeed implies that $\mu$ is inner-regular.  We will also prove that there is a Borel measure $\mu$ on $M$ such that $\int d(m,0) d |\mu|(m)< \infty$ but $\BL\mu \notin \LF(M)$ if and only if the weight of $M$ is greater than the least real-valued measurable cardinal. Therefore the existence of a metric space which admits a counterexample to Problem 2 in \cite{Aliaga_integral} cannot be proven in ZFC.
    
    The results in this article also provide a partial solution to the problem  weather every sequentially normal functional $\phi \in \LF(M)^{**}$ is normal, posed in \cite{Aliaga}.

\section{Preliminaries}

\subsection{Notation}
    This subsection introduces some basic notation.
    \begin{itemize}
        \item We will denote by $w(M)$ the weight of a topological space $M$, which for metrizable spaces coincides with its density character.
        \item We will denote the distance function on a metric space by $d$.
        \item  We will denote the fact that two metric spaces $M$ and $N$ are isometric by $M \isometric N$. 
        \item The measures in this article are assumed to be possibly infinite signed measures unless otherwise stated (see for example \cite{Folland} for the necessary theorems).
        \item If $\mu$ is a measure and $f$ is a $\mu$-integrable function we will denote by $fd\mu$ the measure given by $fd\mu(E)=\int_Efd\mu$.
    \end{itemize}
    
\subsection{Lipschitz free Spaces}
 
This subsection lists the basic definitions and theorems about Lipschitz free spaces. All results in this subsection can be found on \cite{Weaver}. We begin by generalizing the concept of the origin of a Banach space to an arbitrary metric space:

\begin{Def}
    \label{p_metric}
    A \emph{pointed metric space} is a pair $(M,0)$, where $M$ is a metric space and $0 \in M$. The point $0$ is named the \emph{base point of $M$}. We will also denote $\rho=d(\cdot,0)$.
\end{Def}

\begin{obs}
    The base point of every pointed metric space will be denoted by $0$ unless otherwise stated.
\end{obs}

We define the morphisms between pointed metric spaces $M$ and $N$ as the set $\Lip_0(M,N)$ defined in the following way:

\begin{Def}
    \label{Lip_const}
    Let $M$ and $N$ be a pointed metric spaces and let $f \in N^M$. The \emph{Lipschitz constant of $f$} is defined as:\[
        \norm{f}_L=\sup_{\substack{x,y \in M\\ x \neq y}}\frac{d_N(f(x),f(y))}{d_M(x,y)},
    \]
    if $M$ has at least two points, and $0$ otherwise.

    Define $\Lip_0(M,N)$ to be the set of the functions $f \in N^M$ such that $f(0)=0$ and $\norm{f}_L<\infty$. If $N=\R$ we will write $\Lip_0(M)$ instead of $\Lip_0(M,\R)$.
\end{Def}

The set of morphisms between a pointed metric space and a Banach space is in fact a Banach space.

\begin{prop}
    \label{Lip0}
    Let $M$ be a pointed metric space and let $X$ be a Banach space. Then $\Lip_0(M,X)$ is a Banach space with norm $\norm{\cdot}_L$. 
\end{prop}

The Lipschitz free space over $M$ is the canonical predual of $\Lip_0(M)$. One of the ways of constructing it is to first embed $M$ into $\Lip_0(M)$ in the following way:

\begin{prop}
    \label{delta_m}
    Let $M$ be a pointed metric space. We can define a map $\delta: M \to \Lip_0(M)^*$ by setting $\hprod{f,\delta(m)}=f(m)$ for all $m \in M$. Then $\delta$ is an isometric embedding and $\delta[M \setminus \{0\}] $ is linearly independent.
\end{prop}

We now need to simply take the closed span of $\delta[M]$:

\begin{Def}
    \label{LF}
    Let $M$ be a pointed metric space. The \emph{Lipschitz free space over $M$} is the Banach space defined by $\LF(M)=\cspan(\delta[M])$. 
\end{Def}
The notation $\LF(M)$ does not reference the base point of $M$ because the construction of $\LF(M)$ with different choices of base points leads to isometric spaces

We will use the following properties of Lipschitz free spaces:

\begin{thm}
    \label{LF_propr}
    Let $M$ be a pointed metric space. Then $\LF(M)$ has the following properties:
    
    \begin{enumerate}[label=(\alph*)]
        \item $\LF(M)^*\isometric \Lip_0(M)$. We will identify $\LF(M)^*$ with $\Lip_0(M)$.
        \item  The weak$^*$ topology it induces coincides with the topology of pointwise convergence on bounded sets.
        \item  For every Banach space $X$ and every $f\in \Lip_0(M,X)$ there is a map $F \in \BL(\LF(M),X)$ such that $f=F\circ\delta$, that is the following diagram commutes:\[
            \begin{tikzcd}
            M \arrow[rr, "f"] \arrow[dd, "\delta"'] &  & X \\
                                                    &  &   \\
            \LF(M) \arrow[rruu, "F"']               &  &  
            \end{tikzcd}
        \]
        and moreover $\norm{F}=\norm{f}_L$.
        
        If $X = \R$ we will write $\hprod{x,f}$ instead of $\hprod{x,F}$, where $x \in \LF(M)$.
    \end{enumerate}
\end{thm}

\subsection{Real-valued measurable cardinals}

We begin this section with a procedure to extend a measure defined on a subset of a larger space by adding null measure sets:

\begin{Def}
    \label{measure_extension_superspace}
    Let $(M,\Sigma)$ be a measure space, let $N \subset M$ and let $\mu$ be a measure on the measurable space \[
        (N,\{E \cap N \;:\: E \in \Sigma \}).
    \]
    Define $\ext{\mu}{M} (E)=\mu(E \cap N)$ for all $E \in \Sigma$, notice that $\ext{\mu}{M} $ is a measure on $(M,\Sigma)$. If $N \in \Sigma$ we also have that $\ext{\mu}{M}$ is concentrated on $N$.
\end{Def}

\begin{obs}
    $\ext{\mu}{M}$ is an extension in the set-function sense if and only if $N \in \Sigma$. If $N \notin \Sigma$ then $\ext{\mu}{M}$ is an extension in a looser sense.
\end{obs}

We can now introduce one of the main subjects of this article:

\begin{Def}
    \label{nontrivial_measure_cardinal}
    Let $\kappa$ an be infinite cardinal. A \emph{nontrivial measure\footnote{We will use the terminology in \cite{Jech} for measures on cardinals.} on $\kappa$} is a probability measure $\mu$ on $(\kappa,\mathcal{P}(\kappa))$ such that $\mu(\{\alpha\})=0$ for all $\alpha < \kappa$.
\end{Def}

By the extension procedure in \Cref{measure_extension_superspace}, if $\kappa$ is the least cardinal which can have a nontrivial measure on it, then a cardinal $\alpha$ can have a nontrivial measure on it if and only if $\alpha \geq \kappa$. The cardinal $\kappa$ also must have a stronger measurability property:

\begin{Def}
    \label{real_measure_cardinal}
    Let $\kappa$ and $\alpha$ be infinite cardinals. 
    \begin{enumerate}[label=(\alph*)]
        \item An \emph{$\alpha$-additive probability measure on $\kappa$} is a probability measure $\mu$ on $(\kappa,\mathcal{P}(\kappa))$ such that if $\beta < \alpha$ and $(E_\gamma)_{\gamma < \beta }$ is family of pairwise disjoint subsets of $\kappa$ then:\[
            \mu\left( \bigcup\limits_{\gamma < \beta}E_\gamma \right)=\sum\limits_{\gamma < \beta} \mu(E_\gamma).
        \]
        \item $\kappa$ is \emph{real-valued measurable} if there is a $\kappa$-additive nontrivial probability measure on $\kappa$. 
    \end{enumerate}
\end{Def}

\begin{thm}[\cite{Jech} Corollary 10.7]
    \label{nontrivial_measure_cardinal_rv_measurable}
    Let $\kappa$ be the smallest cardinal such that there is a nontrivial measure on $\kappa$. Then $\kappa$ is real-valued measurable.
\end{thm}

From real-valued measurability we get the following property:

\begin{thm}[\cite{Jech} Corollary 10.15]
    \label{rv_measurable_cardinal_weakly_inac}
    Let $\kappa$ be real-valued measurable. Then:
    \begin{enumerate}[label=(\alph*)]
        \item $\kappa$ is uncountable.
        \item There is no cardinal $\alpha$ such that $\kappa$ is the smallest cardinal greater than $\alpha$.
        \item $\kappa$ is \emph{regular}. That is, if $(A_i)_{i \in I}$ is a family of sets such that $|I| < \kappa$ and $|A_i| < \kappa$ for all $i \in I$, then $\left| \bigcup\limits_{i \in I}A_i \right| <\kappa$.
    \end{enumerate}
\end{thm}

As noted on \cite{Jech} page 33, the existence of a cardinal with properties (a), (b) and (c) cannot be proven in ZFC, therefore the existence of a cardinal that admits a nontrivial measure cannot be proven in ZFC.

We can now prove a theorem extracted from the proof of \cite{Bogachev} Proposition 7.2.10 which shows that if the weight of a metric space $M$ allows a nontrivial measure then we can construct an analogous Borel measure on $M$.
\begin{cor}
    \label{measurable_metric}
    Let $M$ be a metric space such that there is a nontrivial measure on the weight of $M$. Then there is a uniformly discrete $M' \subset M$ and a probability measure $\mu$ on $(M',\mathcal{P}(M'))$ such that $\mu({\{m\}})=0$ for all $m \in M'$.
\begin{proof}
    Let $\kappa$ be the smallest cardinal such that there is a nontrivial measure on $\kappa$. For all $n \in \N$ we can use Zorn's lemma to construct a maximal subset $M_n$ such that $d(x,y) \geq n^{-1}$ for all distinct $x,y \in M_n$, that is $d(x,M_n) < n^{-1}$ for all $x \in M$.  Define $M_d =\bigcup\limits_{n \in \N}M_n$; we clearly have that $M_d$ is dense. The set of balls centered on elements of $M_d$ with rational radius is a topological basis of $M$ with cardinality $ |M_d|$. That is $|M_d| \geq w(M) \geq \kappa$, and therefore by \Cref{rv_measurable_cardinal_weakly_inac}(b), $|M_{n_0}| \geq \kappa$, for some $n_0 \in \N$. 
    
    Take $M' \subset M_{n_0}$ with $|M'|=\kappa$ and let $f:M' \to \kappa$ be bijective. We can define the desired measure on $(M',\mathcal{P}(M'))$ by taking $\mu'$ a nontrivial measure on $\kappa$ and setting $\mu(E)=\mu'(f[E])$ for all $E \subset M'$.
    
\end{proof}
\end{cor}
\begin{obs}
    Because $M'$ is uniformly discrete $\mathcal{P}(M') \subset \mathbb{B}_M$, and thus $\ext{\mu}{M}$ is a genuine extension of $\mu$ to a Borel measure on $M$.
\end{obs}

We conclude this section with a theorem that implies that finite (or more generally $\sigma$-finite) measures on a metric space whose weight is smaller than the least real-valued measurable cardinal have strong regularity properties.

\begin{thm}[\cite{Bogachev} Proposition 7.2.10]
    \label{Radon_nomeasurable_equiv}
    Let $M$ be a metric space. Then the following are equivalent:
    \begin{enumerate}[label=(\alph*)]
        \item There is no nontrivial measure on the weight of $M$.
        \item Every finite Borel measure on $M$ is concentrated on a separable subset\footnote{The condition in the original theorem is that every Borel measure $\mu$ on $M$ is concentrated on a specific separable subset called the \emph{support of $\mu$}.} of $M$.
    \end{enumerate}
\end{thm}

\begin{obs}
    By taking a countable partition of $M$ into sets of finite measure we can generalize the condition (b) in \Cref{Radon_nomeasurable_equiv} to $\sigma$-finite measures.
\end{obs}

By \cite{Bogachev} Chapter 7 a finite measure on a complete metric space is Radon if and only if it is concentrated on a separable subset. We can easily generalize the theorems in this chapter to prove that a $\sigma$-finite measure on a complete metric space is inner-regular if and only if it is concentrated on a separable subset.

\section{Measures and Functionals on $\Lip_0(M)$}

In the construction of the space $\LF(M)$ the functional $\delta(m)$ can be interpreted as integration with respect to the Dirac measure $\delta_m$, therefore the space $\LF(M)$ can be regarded as the completion of the space of finitely supported Borel measures on $M \setminus \{0\}$. The following definition will allow us to induce elements of $\LF(M)$ via Borel measures besides finitely supported ones:

\begin{Def}
    \label{Lmu}
    Let $M$ be a pointed metric space and $\mu$ be a Borel measure on $M$ such that $\Lip_0(M) \subset L^1(\mu)$. Define $\BL\mu \in \Lip_0(M)^\#$ (the algebraic dual of $\Lip_0(M)$) by:\[
        \BL\mu(f)=\int f d\mu.
    \]
\end{Def}

Let $M$ be a pointed metric space and let $\overline{M}$ denote its completion. Then we can define a weak$^*$-weak$^*$ continuous isometric isomorphism $T \in \BL(\Lip_0(\overline{M}),\Lip_0(M))$ by setting $T(f)=\rest{f}{M}$ for all $\Lip_0(\overline{M})$. In particular we have that $\LF(\overline{M})\isometric\LF(M)$. For this reason it is common in the theory of Lipschitz free spaces to assume $M$ is complete, as done, for example, in \cite{Aliaga_integral}. 

We will prove a proposition that allows us to apply the theorems in \cite{Aliaga_integral} used in this chapter to non complete metric spaces:

\begin{prop}
    Let $M$ be a pointed metric space and let $\mu$ be a Borel measure on $M$. Then:
    \begin{enumerate}[label=(\alph*)]
        \item 
        For every $f \in \Lip_0(\overline{M})$ we have that:\[
            \int \rest{f}{M} d \mu =\int f d( \ext{\mu}{\overline{M}}).
        \]
        \item 
        $\BL\mu \in \LF(M)$ if and only if $\BL( \ext{\mu}{\overline{M}}) \in \LF(\overline{M})$.
    \end{enumerate}
\begin{proof}
    \begin{enumerate}[label=(\alph*)]
        \item 
            By proceeding with $\mu^+,\mu^-$ and $f^+,f^-$, we can suppose without loss of generality that $\mu$ and $f$ are positive. 

            First, given $E \in \mathbb{B}_{\overline{M}}$ we have that:
            \[
                \int \rest{\charc_E}{M} d \mu= \mu(E \cap M)= \ext{\mu}{\overline{M}}(E)=\int \charc_E  d( \ext{\mu}{\overline{M}}).
            \]
            Now, denote by $S$ the set of the simple functions $s \in \R^{\overline{M}}$ such that $s \leq f$ $\ext{\mu}{\overline{M}}$-almost-everywhere and by $S_M$ the set of the simple functions $s \in \R^{{M}}$ such that $s \leq \rest{f}{M}$ $\mu$-almost-everywhere. We can easily prove that:\[
                S_M=\{\rest{s}{M}: \; s \in S\},
            \]
            therefore:\[
                \int \rest{f}{M} d \mu =\sup_{s \in S}\int \rest{s}{M} d\mu =\sup_{s \in S_M}\int s d( \ext{\mu}{\overline{M}}) =\int f d( \ext{\mu}{\overline{M}}).
            \]
    \item 
        By (a) we have that $\Lip_0(M) \subset L^1(\mu)$ if and only if $\Lip_0(\overline{M}) \subset L^1(\ext{\mu}{\overline{M}})$ and that $\BL\mu$ is $\omega^*$-continuous if and only if $ \BL( \ext{\mu}{\overline{M}})$ is $\omega^*$-continuous. Therefore $\BL\mu \in \LF(M)$ if and only if $\BL( \ext{\mu}{\overline{M}}) \in \LF(\overline{M})$.
    \end{enumerate}
        
\end{proof}
\end{prop}

We can easily see that $\BL\mu$ is well defined if and only if $\mu$ induces a linear functional:

\begin{prop}[\cite{Aliaga_integral} Proposition 4.3]
    \label{Lmu_equiv}
    Let $M$ be a pointed metric space and let $\mu$ be a Borel measure on $M$. Then the following conditions are equivalent:
    \begin{enumerate}[label=(\alph*)]
        \item $\Lip_0(M) \subset L^1(\mu)$.
        \item $\int \rho d|\mu|<\infty$.
        \item $\BL\mu \in \LF(M)^ {**}$.        
    \end{enumerate}

\end{prop}

If $\BL\mu \in \LF(M)^{**}$ we have the following consequences:

\begin{cor}[\cite{Aliaga_integral} Proposition 4.3]
    \label{Lmu_propr}
    Let $M$ be a pointed metric space and let $\mu$ be a Borel measure on $M$ such that $|\mu|(\{0\}) < \infty$ and $\BL\mu \in \LF(M)^{**}$. Then:
    \begin{enumerate}[label=(\alph*)]
        \item $|\mu|(K)<\infty$ for all closed $K \subset M \setminus \{0\}$.
        \item $\mu$ is $\sigma-$finite.
    \end{enumerate}
\end{cor}

We can now show an example\footnote{This is also basically the standard example of an incomplete metric space $M$ such that there is a nonzero $m \in \LF(M)$ having empty support. See \cite[Proposition 3.3.10]{Aliaga}.} of a pointed metric space $M$ and an element of $\LF(M)$ that is induced by a Borel measure on $\overline{M}$ but is not induced by a Borel measure on $M$:

\begin{ex}
    \label{meas_on_complt}
    Let $M=\{0\} \cup (1,2]$. By a slight abuse of notation we have that $\BL\delta_1=\delta(1) \in \LF(M)$. But if $\mu$ is a Borel measure on $M$ such that $\BL\mu \in \LF(M)$ then $\BL\mu \neq \delta(1)$.

    \begin{proof}
        Suppose by contradiction that $\mu$ is a Borel measure on $M$ such that $\BL(\mu)=\delta(1)$.
        
        If we define $\mu_0$ by $d\mu_0=\charc_{M \setminus \{0\}}d\mu$ then $\BL\mu=\BL\mu_0$. Therefore we may assume without loss of generality that $\mu$ is concentrated on $(1,2]$. First notice that $\mu$ is finite  by \Cref{Lmu_propr}(a). Now fix a closed $F \subset (1,2]$ and denote $k=d(1,F)^{-1}$. We can define $f_n \in \Lip_0(M)$ by:\[
        f_n(m)= \max\Big(0, \big(1-k2^{n} d( m,F) \big) \Big)
        \]
        for all $m \in M$. Therefore, by the Lebesgue dominated convergence theorem and the Radon-Nikodym theorem, we have that:\[
            \mu(F)=\int \charc_F d \mu=\lim_{n \to \infty}\int f_n d \mu=0
        \]
        But, $\mu$ is finite and therefore regular (see for example \cite[Theorem 7.1.7]{Bogachev}). From which we conclude that $\mu=0$, a contradiction.
    \end{proof}
\end{ex}

We can now show a sufficient and necessary condition for a Borel measure to induce an element of $\LF(M)$:

\begin{thm}
    \label{Lmu_separability}
    Let $M$ be a pointed metric space and let $\mu$ be a Borel measure on $M$. Then $\BL\mu \in \LF(M)$ if and only if $\int \rho d|\mu|<\infty$ and $\mu$ is concentrated on a separable subset of $M$.
\begin{proof}
    By \Cref{Lmu_equiv} we only need to show that if $M$ is a pointed metric space and $\mu$ is a Borel measure on $M$ such that $\int \rho d|\mu|<\infty$, then $\BL\mu \in \LF(M)$ if and only if $\mu$ is concentrated on a separable subset of $M$.
    \begin{enumerate}
        \item [$\implies)$]
        Notice that if we define $\mu_0$ by $d\mu_0=\charc_{M \setminus \{0\}}d\mu$ then $\BL\mu=\BL\mu_0$. Therefore we may assume without loss of generality that $\mu(\{0\})=0$.
        
        By the hypothesis that $\BL\mu \in \LF(M)$, we can write $\BL\mu$ in the form:\[
            \BL\mu=\lim_{n_1 \to \infty}\sum_{n_2 \in \N}a_{ (n_1,n_2)}\delta(m_{(n_1,n_2)}),
        \]    
        where $a_n \in \R$ and $m_n\in M$ for all $n \in \N^2$ and for each $n_1 \in \N$ the set of the $n_2 \in \N$ for which $a_{(n_1,n_2)} \neq 0$ is finite.
        
        Fix $i \in \N$ and define:
        \begin{align*}
            M_i &=\rho^{-1}[[2^{-i},\infty)],\\
            M'&=\overline{\{m_n: n \in \N^2\}}  \cup \{0\}. 
        \end{align*} 
        For each $k \in \N$ define ${f_{i,k}^u, f_{i,k}^l, f_{i,k} \in \Lip_0(M)}$ by:
        \begin{align*}
            f_{i,k}^u(m)&=\max\Big(0, \big(1-k2^{i+1} d( m,M_{i}) \big) \Big)\\
            f_{k}^l(m)&=kd(m,M')\\
            f_{i,k}  (m)&=\min(f_{i,k}^u(m), f_{k}^l(m) )\\
        \end{align*}
        for all $m \in M$. Notice that $\rest{f_{i,k}}{M'}=0$ for all $k \in \N$. Also, letting $k \to \infty$, we have that $ f_{i,k}$ converges pointwise to $\charc_{M_i \setminus M'}$, and for all $i \in \N$ we have that $ |f_{i,k}| \leq |f_{i,k}^u| \leq \charc_{M_{i+1}} \in L^1(\mu)$ by \Cref{Lmu_propr}, therefore by the Lebesgue Dominated Convergence Theorem and the Radon-Nikodym Theorem:
        \begin{equation*}            
        \begin{aligned}
            \mu(M_i \setminus M')&= \int \charc_{M_i \setminus M'}d\mu = \lim_{k \to \infty} \hprod{\BL\mu,f_{i,k}}\\
                               &= \lim_{k \to \infty} \hprod{\lim_{n_1 \to \infty}\sum_{n_2 \in \N}a_{ (n_1,n_2)}\delta(m_{(n_1,n_2)}),f_{i,k}}\\
                               &= \lim_{k \to \infty} \lim_{n_1 \to \infty}\sum_{n_2 \in \N}a_{ (n_1,n_2)}f_{i,k}(m_{(n_1,n_2)})\\
                               &=0.
        \end{aligned}
        \end{equation*}
        We have that $(M_i \setminus M')_{i \in \N}$ increases to $M \setminus M'$, therefore $\mu(M \setminus M')=0$ by the continuity of the measure. We will conclude by proving that $|\mu|(M \setminus M')=0$.
        
        Fix $A^+,A^-$ a Hahn decomposition associated to $\mu$ and  $\epsilon>0$, also denote  $B^\pm=A^\pm\cap (M \setminus M')$. Notice that $\mu(B^+)=\mu(B^-)$ and therefore ${|\mu|(M \setminus M')< \infty}$.
        
        Define $\nu$ by $d\nu=\charc_{M \setminus M'}d\mu$; we have that $\nu$ is finite, therefore because $\nu$ is also regular (see for example \cite[Theorem 7.1.7]{Bogachev}) we can take closed sets $F^\pm \subset B^\pm$ and open sets $V^\pm \supset B^\pm$ such that $|\nu|(V^\pm\setminus F^\pm)<\epsilon$. Also, because $B^\pm \subset M \setminus M'$, which is open, we can assume $V^\pm \subset M \setminus M'$, and thus: \[
            |\mu|(V^\pm\setminus F^\pm)=|\nu|(V^\pm\setminus F^\pm)<\epsilon.
        \]
        We have that $F^+$, $F^-$ and $K = M \setminus(V^+ \cup V^-) \cup \{0\}$ are closed and pairwise disjoint, therefore we can take $g \in \Lip_0(M)$ such that ${F^+ \subset g^{-1}[1]}$, ${F^- \subset g^{-1}[-1]}$, $K \subset g^{-1}[0]$, and $|g(m)| \leq 1$ for all $m \in M$. Thus we have that:
        \begin{align*}
            |\mu|(M \setminus M')&=|\mu|(B^+) + |\mu|(B^-) \\
                                 &\leq 2\epsilon +  |\mu|(F^+) + |\mu|(F^-) \\
                                 &=2 \epsilon + \int_{ F^+ \cup F^- \cup K } g  d\mu \\
                                 &\leq 4 \epsilon +\int g d \mu.
        \end{align*}
        Finally, for each $k \in \N$ define $h_k \in \Lip_0(M)$ by:\[
            h_{k}(m)=\min(1,kd(m,M'))
        \]
        for each $m \in M$. Notice that $\rest{h_{k}}{M'}=0$ for all $k \in \N$ and that $h_k$ increases pointwise to $\charc_{M \setminus M'}$. Therefore, because $V^\pm \subset M \setminus M'$ and $\rest{g}{M \setminus( V^+ \cup V^-)}=0$, we have that $gh_k$ converges pointwise to $g$ and that $|gh_k| \leq \charc_{M \setminus M'} \in L^1(\mu)$ for all $k \in \N$. That is, by the Lebesgue Dominated Convergence Theorem and the Radon-Nikodym Theorem:\[
            \int g d\mu=\lim_{k \to \infty}\int h_{k}g d\mu = \lim_{k \to \infty} \hprod{\BL\mu,h_{k}g } =0, 
        \]
        thus, we get the desired result by letting $\epsilon \to 0$.
        \item [$\impliedby)$]
        By hypothesis we can take a  $E \in \B_M$ such that ${|\mu|(E)=0}$ and $M \setminus E$ is separable. Therefore the result follows from \cite{Aliaga_integral} Proposition 4.4.
    \end{enumerate}
\end{proof}
\end{thm}

 To fully answer Problem 2 of \cite{Aliaga_integral} it only remains to be proved that a metric space $M$ and a Borel measure $\mu$ on $M$ failing the conditions in \Cref{Lmu_separability} indeed exist. But, as we will show in the next theorem, the existence of such a measure is equivalent to the existence of a real-valued measurable cardinal. That is, assuming that ZFC is consistent, the existence of $\mu$ cannot be proved and possibly cannot be refuted.

\begin{cor}
    \label{Lmu_M_nonmeasurable_w}
    Let $M$ be a pointed metric space. Then the following are equivalent:
    \begin{enumerate}[label=(\alph*)]
        \item There is no nontrivial measure on the weight of $M$.
        \item Every $\sigma$-finite Borel measure on $M$ is concentrated on a separable subset of $M$.
        \item For every  Borel measure $\mu$ on $M$ such that $\int \rho d|\mu|<\infty$ we have that $\BL\mu \in \LF(M)$.
    \end{enumerate}
\begin{proof}
    We know that (a) and (b) are equivalent by \Cref{Radon_nomeasurable_equiv}, also, (b) implies (c) by \Cref{Lmu_separability}. It only remains to be proven that (c) implies (a), which we will do by contraposition.
    
    By \Cref{measurable_metric} we can take  a $M' \subset M$ such that there is a $r>0$ for which $d(x,y) \geq 3r$ for all $x,y \in M'$ and there is a probability measure $\mu'$ on $(M',\mathcal{P}(M'))$ for which $\mu(\{x\})=0$ for all $x \in M'$. We can assume without loss of generality that $d(0,M') > r$. Define $f \in \Lip_0(M)$ by: \[
        f(m)=r-  \min(d(m,M'),r),
    \]
    for all $m \in M$. And, given $N \in \mathcal{P}_{\fin}(M)$, define $f_N \in \Lip_0(M)$ by:\[
        f_N(m)=r- \min(d(m,N),r) ,
    \]
    for all $m \in M$. Clearly $f_N$ converges pointwise and boundedly to $f$, therefore $f_N \xrightarrow{\omega^*}f$.
    
    Define $\mu$ by $d\mu=\frac{1}{\rho}d(\ext{\mu'}{M}) $ and for each $n \in \Z$ define $C_n=\rho^{-1}[[2^n,2^{n+1} )] $. Clearly $\int \rho d|\mu|<\infty$ and $\ext{\mu'}{M}(C_i)>0 $ for some $ i \in \Z$, therefore:\[
        \lim\limits_{N \in \mathcal{P}_{\fin}(M)}\int f_N d\mu=0 < \int_{C_i}r2^{ -(i+1)}\charc_{M'}d (\ext{\mu'}{M}) \leq \int \frac{f}{\rho} d (\ext{\mu'}{M}) = \int f d\mu.
        \]
    That is, $\BL\mu$ is not $\omega^*$-continuous, so $\BL\mu \notin \LF(M)$.
\end{proof}
\end{cor}

\section{Normal Functionals on $\Lip_0(M)$}

The mapping $\BL$ is related to the concept of normal functionals on $\Lip_0(M)$:

\begin{Def}
    \label{normal_func}
    Let $M$ be a pointed metric space. A functional $\phi \in \LF(M)^{**}$ is \emph{normal} if, given a norm bounded net $(f_i)_{i \in I}$ on $\Lip_0(M)$ that converges pointwise and monotonically to some $f \in \Lip_0(M)$, we have that $\hprod{f_i,\phi} \to \hprod{f,\phi}$.
    
    The functional $\phi$ is \emph{sequentially normal} if given a norm bounded sequence $(f_i)_{i \in \N}$ on $\Lip_0(M)$ that converges pointwise and monotonically to some $f \in \Lip_0(M)$ we have that $\hprod{f_i,\phi} \to \hprod{f,\phi}$.
    
\end{Def}

The concept of normality can be used as a characterization of weak$^*$ continuity:

\begin{thm}[\cite{Aliaga_normal} Theorem 2] 
    \label{normal_func_is_ws_cont}
    Let $M$ be a pointed metric space and let $\phi \in \LF(M)^{**}$. Then $\phi$ is normal if and only if it is weak$^*$ continuous.
\end{thm}

\Cref{Lmu_M_nonmeasurable_w} provides a partial answer to a question about normality of functionals on $\Lip_0(M)$:  

{\cite{Aliaga} Problem 3.2}

\begin{quote}
    Is sequential normality equivalent to normality, and hence to weak$^*$ continuity?
\end{quote}

Namely, if $\phi=\BL\mu$ for some Borel measure $\BL\mu$ it is sequentially normal by the Lebesgue Dominated Convergence Theorem and the Radon-Nikodym Theorem, but if the weight of $M$ is greater than the least real-valued measurable cardinal then there is a Borel measure $\mu$ on $M$ such that $\BL\mu \in \LF(M)^{**} \setminus \LF(M)$. That is, given the existence of a real-valued measurable cardinal we can construct a counterexample to \cite{Aliaga} Problem 3.2.

An immediate follow up question is whether the existence of a real-valued measurable cardinal is necessary to construct such a counterexample, or more generally:
\begin{prob}
    Can the existence of a metric space $M$ such that there is a sequentially normal functional $\phi \in \LF(M)^{**} \setminus \LF(M)$ be proven in ZFC?
\end{prob}

%%%%%%%%%%% To ease editing, use normal size for the references:

\end{document}